\documentclass[reqno,a4paper,12pt]{amsart}

\ifdefined\SMART
\usepackage[paperwidth=12cm,paperheight=24cm,left=5mm,right=5mm,top=10mm,bottom=5mm]{geometry}
\else
\usepackage[margin=1in]{geometry}
\fi

\usepackage[T1]{fontenc}
\usepackage[utf8]{inputenc}
\usepackage{amsmath,amssymb,mathtools}
\usepackage{enumerate}
\usepackage[sort,nocompress]{cite}
\usepackage[hidelinks]{hyperref}
\numberwithin{equation}{section}
\mathtoolsset{showonlyrefs}
\newtheorem{definition}{Definition}[section]
\newtheorem{theorem}[definition]{Theorem}
\newtheorem{lemma}[definition]{Lemma}
\newtheorem{proposition}[definition]{Proposition}
\newtheorem{corollary}[definition]{Corollary}
\newcommand{\N}{\mathbb N}
\newcommand{\R}{\mathbb R}
\newcommand{\C}{\mathbb C}
\newcommand{\K}{\mathbb K}
\newcommand{\T}{\mathbb T}
\newcommand{\A}{\mathcal A}
\newcommand{\Pro}{\mathcal P}
\newcommand{\norm}[1]{\left\lVert#1\right\rVert}

\allowdisplaybreaks
\usepackage{xcolor}

\title[Restricted forms of Fekete's lemma and reflexivity]{Restricted forms of Fekete's lemma and reflexivity in Banach spaces}
\hypersetup{pdftitle={Restricted forms of Fekete's lemma and reflexivity in Banach spaces}}

\author[V.~Blakaj]{Vjosa Blakaj}
\address[Vjosa Blakaj]{Department of Mathematical Sciences, University of Copenhagen, Universitetsparken 5, 2100 Copenhagen Ø, Denmark}
\email{vb@math.ku.dk}

\author[S.~Grepstad]{Sigrid Grepstad}
\address[Sigrid Grepstad]{Department of Mathematical Sciences, Norwegian University of Science and Technology (NTNU), NO-7491 Trondheim, Norway}
\email{sigrid.grepstad@ntnu.no}

\author[M.~Kolountzakis]{Mihail N. Kolountzakis}
\address[Mihalis Kolountzakis]{{Department of Mathematics and Applied Mathematics}, University of Crete, Voutes Campus, 70013 Heraklion, Greece,\newline and\newline {Institute of Computer Science, Foundation of Research and Technology Hellas}, N. Plastira 100, Vassilika Vouton, 70013 Heraklion, Greece}
\email{kolount@uoc.gr}

\author[L.~Liehr]{Lukas Liehr}
\address[Lukas Liehr]{Department of Mathematics, Bar-Ilan University,
Ramat-Gan 5290002, Israel}
\email{lukas.liehr@biu.ac.il}

\date{}

\begin{document}
\begin{abstract}In this work, we study restricted forms of Fekete’s lemma for vector-valued sequences, in which the defining inequalities are imposed only on a prescribed set of index pairs. We identify a condition on weighted graph distances that is necessary and sufficient for the lemma to hold in the Euclidean plane and sufficient in every Hilbert space. Using a renorming theorem of Odell and Schlumprecht, we show that a separable Banach space of real dimension at least two admits an equivalent norm with the restricted Fekete property if and only if it is reflexive and the prescribed set satisfies the graph distance condition. Moreover, a single norm can be chosen for all such sets.
\end{abstract}



\maketitle

\section{Introduction and main results}\label{sec:intro}

Fekete's lemma \cite{fekete1923} states that if a real sequence $(a_n)_{n\in\N}$ is subadditive, that is, for all $n,m\in\N$ we have $a_{n+m}\leq a_n+a_m$, then
$$
\lim_{n\to\infty}\frac{a_n}{n}=\inf_{n\in\N}\frac{a_n}{n}\in[-\infty,\infty).
$$
A vector-valued extension of Fekete's lemma was proved by Kulikov and Shao \cite{kulikov2024feketeslemmabanachspaces}: if $X$ is a uniformly convex Banach space and $(v_n)_{n\in\N}\subset X$ satisfies
\begin{equation}\label{eq:fekete}
\norm{v_{n+m}}\leq\norm{v_n+v_m}, \quad n,m\in\N,
\end{equation}
then $v_n/n$ converges in norm to an element of $X$.

The scalar and vector statements behave differently when the inequalities are restricted. De Bruijn and Erd\H{o}s \cite[Theorem 22]{deBruijnErdos1952} proved that the scalar version of Fekete's lemma still holds if subadditivity is assumed only for $n\leq m\leq2n$. In \cite{deBruijnErdos1952} the authors also claimed that the number 2 in the inequality $n \le m \le 2n$ cannot be reduced, yet F\"uredi and Ruzsa \cite{furedi2020nearly} corrected this claim by showing that the number 2 can be replaced by any constant larger than 1.
Kulikov and Shao showed in \cite[Section 3]{kulikov2024feketeslemmabanachspaces} that the corresponding statement fails already in the Euclidean plane. Moreover, they posed the question whether the upper bound $2n$ can be replaced by $g(n)$ for some $g:\N\to\N$ so that convergence holds.

In what follows, we study which pairs $(n,m)$ need to be compared in \eqref{eq:fekete} to force convergence of $v_n/n$. In addition, we ask which Banach spaces admit an equivalent norm for which these restricted inequalities imply convergence. Our results separate these questions into a condition on the comparison set and a geometric condition on the underlying Banach space. The former is expressed using weighted graphs, whereas the latter is expressed through reflexivity of the underlying space.

\subsection{Restricted Fekete property}

Throughout, Banach spaces are over the real or complex numbers, and convergence of vectors means convergence in norm unless stated otherwise. Since \eqref{eq:fekete} is symmetric in $n$ and $m$, we consider sets of the form
$$
\A\subseteq\Delta, \qquad \Delta=\{(n,m)\in\N^2:n\leq m\}.
$$
For a fixed norm $\norm{\cdot}$ on $X$, a sequence $(v_n)_{n\in\N}\subset X$ is called \emph{$\A$-admissible} if
$$
\norm{v_{n+m}}\leq\norm{v_n+v_m}, \quad (n,m)\in\A.
$$
We say that the norm $\norm{\cdot}$ has the \emph{$\A$-Fekete property} if, for every $\A$-admissible sequence, there exists $v\in X$ such that
$$
\lim_{n\to\infty}\norm{v_n/n-v}=0.
$$
If $\A=\Delta$ we talk about the \emph{unrestricted Fekete property}.

The \emph{$\A$-Fekete property} depends on the norm used in the inequalities. Recall that for a Banach space $(X,\norm{\cdot}_X)$, a norm $q$ is equivalent to $\norm{\cdot}_X$ if there are constants $0<c\leq C<\infty$ such that
$$
c\norm{x}_X\leq q(x)\leq C\norm{x}_X, \quad x\in X.
$$
Thus, saying that $X$ admits an equivalent $\A$-Fekete norm means that there exists a norm $q$ so that for all sequences $(v_n)_{n \in \mathbb{N}}$ with
$$
q(v_{n+m})\leq q(v_n+v_m), \quad (n,m)\in\A,
$$
one has that $v_n/n$ converges to some element $v \in X$ with respect to the norm $q$.

\subsection{A weighted graph}\label{subsec:graph}

Fix $\A\subseteq\Delta$. For each real sequence $b=(b_n)_{n \in \mathbb{N}}$ satisfying
\begin{equation}\label{eq:scalar-profiles}
0\leq b_n\leq n, \quad n\in\N, \qquad b_{n+m}\leq b_n+b_m, \quad (n,m)\in\A,
\end{equation}
consider the undirected graph $G=(V,E)$ where $V= \N$ and $(n,m) \in E$ if and only if $(\min(n,m),\max(n,m))\in\A$. Given an edge $(n,m) \in E$, we give it the cost
$$
w_b(n,m)=\sqrt{\frac{b_n+b_m-b_{n+m}}{\min(n,m)}}.
$$
Notice that the quantity $b_n+b_m-b_{n+m}$ can be regarded as a subadditivity defect. For $r,s\in\N$, let $d_b(r,s)$ be the infimum of the total edge costs over all finite paths from $r$ to $s$, allowing repeated vertices. We set $d_b(r,r)=0$ and $d_b(r,s)=\infty$ if no such path exists.

\subsection{Defect condition}\label{subsec:defect}

We say that $\A$ satisfies the \emph{defect condition} if, for every $b$ satisfying \eqref{eq:scalar-profiles} we have
\begin{equation}\label{eq:defect-condition}
\lim_{N\to\infty}\sup_{r,s\geq N}d_b(r,s)=0.
\end{equation}
Equivalently, for each such $b$ and each $\varepsilon>0$, there is $N$ such that any two vertices $r,s\geq N$ can be joined by a finite path of total cost less than $\varepsilon$. The path may pass through vertices below $N$, and $N$ may depend on $b$.

The definitions above are motivated by the fact that the subadditivity defect controls the distance between the directions of two vectors in a Hilbert space. To see this, suppose first that $(z_n)_{n \in \mathbb{N}}$ is $\A$-admissible and satisfies $n\leq a_n:=\norm{z_n}\leq2n$. Putting $b_n=a_n-n$ and $u_n=z_n/a_n$ we notice that $b$ satisfies \eqref{eq:scalar-profiles}. Moreover, subtracting $n$ leaves the subadditivity defects unchanged.

For each $(n,m)\in\A$, we have
$$
\begin{aligned}
\norm{u_n-u_m}^2
&=\frac{(a_n+a_m)^2-\norm{z_n+z_m}^2}{a_na_m}\\
&\leq\frac{4(a_n+a_m-a_{n+m})}{\min(a_n,a_m)}\\
&\leq\frac{4(b_n+b_m-b_{n+m})}{\min(n,m)}
=4w_b(n,m)^2.
\end{aligned}
$$
Thus, the defect divided by the smaller index bounds the squared distance between the unit vectors. This explains both the denominator and the square root in the edge cost: $\norm{u_n-u_m}\leq2w_b(n,m)$, on every allowed pair. Summing this estimate along a finite path, and taking the infimum over all such paths, gives
$$
\norm{u_r-u_s}\leq2d_b(r,s), \quad r,s\in\N.
$$
This is the reason for using path distances: they extend the control from allowed pairs to arbitrary indices. Condition \eqref{eq:defect-condition} then makes $(u_n)_{n \in \mathbb{N}}$ a Cauchy sequence.

\subsection{Main results}

Our first result determines when the restricted inequalities characterize the same class of spaces, up to equivalent renorming, as the unrestricted inequalities. The additional requirement is exactly condition \eqref{eq:defect-condition}, which involves only $\A$.

\begin{theorem}\label{thm:transfer}
Let $X$ be a separable Banach space of dimension at least two, and let $\A\subseteq\Delta$. The following are equivalent:
\begin{enumerate}[\rm (i)]
\item $X$ admits an equivalent norm with the $\A$-Fekete property.
\item $\A$ satisfies condition \eqref{eq:defect-condition}, and $X$ admits an equivalent norm with the unrestricted Fekete property.
\end{enumerate}
\end{theorem}

The norms whose existence is asserted in (i) and (ii) may differ. In particular, the theorem concerns the existence of suitable norms. It does not claim that every unrestricted Fekete norm also has the restricted property. For Hilbert spaces we have the following result.

\begin{theorem}\label{thm:test-spaces}
Let $\A\subseteq\Delta$. The following are equivalent:
\begin{enumerate}[\rm (i)]
\item $\A$ satisfies condition \eqref{eq:defect-condition}.
\item The Euclidean norm on $\R^2$ has the $\A$-Fekete property.
\item Every real or complex Hilbert space has the $\A$-Fekete property in its Hilbert norm.
\end{enumerate}
\end{theorem}

Thus, the same condition governs the Euclidean plane and Hilbert spaces of arbitrary dimension, including nonseparable ones. The proof gives a stronger necessity statement: if condition \eqref{eq:defect-condition} fails, then every normed space of real dimension at least two contains an $\A$-admissible sequence for which $v_n/n$ does not converge.\vspace{0.2cm}

We now come to the setting of separable Banach spaces. Recall that $X$ is \emph{reflexive} if the canonical map $J_X:X\to X^{**}$, defined by $(J_Xx)(f)=f(x)$ for $f\in X^*$, is onto. Here $X^*$ is the continuous dual of $X$.

\begin{theorem}\label{thm:reflexivity-characterization}\label{thm:classification}\label{thm:strong-classification}
Let $X$ be a separable Banach space of dimension at least two, and let $\A\subseteq\Delta$. The following are equivalent:
\begin{enumerate}[\rm (i)]
\item $X$ admits an equivalent norm with the $\A$-Fekete property.
\item $X$ is reflexive and $\A$ satisfies condition \eqref{eq:defect-condition}.
\end{enumerate}
Moreover, for each such reflexive space $X$, there is one equivalent norm that has the $\A$-Fekete property for every $\A\subseteq\Delta$ satisfying condition \eqref{eq:defect-condition}.
\end{theorem}

The last assertion makes the order of choices explicit: the norm is chosen from $X$ before the set $\A$ is specified. Taking $\A=\Delta$ gives, in particular, that a separable Banach space admits an equivalent (classical) Fekete norm if and only if it is reflexive. This last statement also holds in real dimensions one. The dimension assumption in the theorem is needed for the necessity of the condition on $\A$.

\section{Sufficiency of the defect condition}\label{sec:hilbert}

In this section, we prove that the defect condition \eqref{eq:defect-condition} is sufficient for the $\A$-Fekete property to hold in Hilbert spaces. Specifically, if $\A$ satisfies condition \eqref{eq:defect-condition}, then $(v_n/n)_{n \in \mathbb{N}}$ converges in norm for every $\A$-admissible sequence $(v_n)_{n \in \mathbb{N}}$ (Proposition \ref{prop:hilbert-sufficiency}).

\subsection{Notation}\label{subsec:notation}

Let $\A \subseteq \Delta$, where $\Delta=\{(n,m)\in\N^2:n\leq m\}$. Recall that a sequence $(v_n)_{n \in \mathbb{N}}$ in a normed space is $\A$-admissible if $\|v_{n+m}\| \leq \|v_n+v_m\|$ for all $(n,m) \in \A$. A sequence $a=(a_n)_{n \in \mathbb{N}}$, where each element $a_n$ is nonnegative, is called \emph{$\A$-subadditive} if $a_{n+m} \leq a_n+a_m$ for all $(n,m) \in \A$. The \emph{defect} of a real sequence $a=(a_n)_{n \in \mathbb{N}}$ at $(n,m) \in \N^2$ is defined by
$$
D_a(n,m) = a_n+a_m-a_{n+m}.
$$
Thus, a nonnegative sequence $a$ is $\A$-subadditive if and only if $D_a(n,m) \geq 0$ for all $(n,m) \in \A$. We denote by $\Pro(\A)$ the set of all real sequences $b=(b_n)_{n \in \mathbb{N}}$ that are $\A$-subadditive with $b_n \leq n$ for all $n \in \N$.\vspace{0.1cm}

Let $b \in \Pro(\A)$ and recall the definition of the graph $G$ from Subsection \ref{subsec:graph}, together with the edge costs $w_b$ and the distances $d_b$ from Subsection \ref{subsec:graph}. 

For $r,s \in \N$, a \emph{path} from $r$ to $s$ is a finite sequence $(n_0,\dots,n_k)$ in $\N$ with $k \in \N$, $n_0=r$ and $n_k=s$ such that
$$
(\min\{n_{j-1},n_j\},\max\{n_{j-1},n_j\}) \in \A, \quad j=1,\dots,k.
$$
The cost of such a path is $\sum_{j=1}^k w_b(n_{j-1},n_j)$. Thus, if $r \neq s$, then $d_b(r,s)$ is the infimum of the costs of all paths from $r$ to $s$. Moreover, $d_b$ satisfies
\begin{equation}\label{eq:pseudometric}
d_b(r,r)=0, \quad d_b(r,s)=d_b(s,r), \quad d_b(r,t) \leq d_b(r,s)+d_b(s,t), \quad r,s,t \in \N.
\end{equation}

If $(n,m) \in \A$, then $D_b(n,m) \geq 0$ and $\min\{n,m\}=n$. Hence, the cost of the edge $(n,m)$ is $w_b(n,m) = (D_b(n,m)/n)^{1/2}$. If $n<m$, then $(n,m)$ is a path from $n$ to $m$, and if $n=m$, then $d_b(n,m)=0$. In both cases we obtain
\begin{equation}\label{eq:one-edge}
d_b(n,m) \leq w_b(n,m) = \sqrt{\frac{D_b(n,m)}{n}}, \quad (n,m) \in \A.
\end{equation}

Notice that the supremum in \eqref{eq:defect-condition} is nonincreasing in $N$. Hence, the limit in \eqref{eq:defect-condition} exists in $[0,\infty]$ for every $b \in \Pro(\A)$, and $\A$ satisfies the defect condition if and only if this limit vanishes for every $b \in \Pro(\A)$.\vspace{0.2cm}

Finally, let $X$ be a normed space over $\K \in \{\R,\C\}$. We denote by $\dim_\R X$ the dimension of $X$ regarded as a real vector space.

\subsection{A linear bound}

We start by showing that the defect condition implies a linear growth bound for $\A$-subadditive sequences.

\begin{lemma}\label{lem:linear-bound}
Let $\A \subseteq \Delta$ satisfy condition \eqref{eq:defect-condition}. Then, there exists $K \in \N$ such that every integer $k>K$ can be written as $k=n+m$ with $(n,m) \in \A$. Moreover, every $\A$-subadditive sequence $a=(a_n)_{n \in \mathbb{N}}$ satisfies
$$
a_n \leq Cn, \quad n \in \N,
$$
where $C = \max_{1 \leq k \leq K} a_k/k$.
\end{lemma}

\begin{proof}
Let $\mathcal N$ be the set of all $k \in \N$ which cannot be written as $k=n+m$ with $(n,m) \in \A$. We first show that $\mathcal N$ is finite.

Suppose by contradiction that $\mathcal N$ is infinite and define $b_n = n\mathbf 1_{\mathcal N}(n)$. For every $(n,m) \in \A$ we have $n+m \notin \mathcal N$, and hence $b_{n+m}=0$. Since $0 \leq b_n \leq n$, it follows that $b \in \Pro(\A)$. Now let $s \in \mathcal N$ and consider an edge of $G$ joining $s$ to some $t \in \N$. Since $b_s=s$ and $b_{s+t}=0$, its defect is $D_b(s,t) = s+b_t \geq s$, while its smaller endpoint satisfies $\min\{s,t\} \leq s$. Hence, every edge incident to $s$ has cost at least one. Since every path from $r$ to $s$ ends with an edge incident to $s$, we obtain $d_b(r,s) \geq 1$ for all distinct $r,s \in \mathcal N$. As $\mathcal N$ is infinite, it follows that $\sup_{r,s \geq N} d_b(r,s) \geq 1$ for every $N \in \N$, which contradicts \eqref{eq:defect-condition}. Hence, $\mathcal N$ is finite.

Since $1 \in \mathcal N$, we may define $K = \max \mathcal N$, and every integer $k>K$ can be written as $k=n+m$ with $(n,m) \in \A$. We prove the bound $a_k \leq Ck$ by induction on $k$, where in the induction step we assume that the bound holds for all indices smaller than $k$. For $k \leq K$ the bound holds by the definition of $C$. If $k>K$, we choose $(n,m) \in \A$ with $k=n+m$. Since $n,m<k$, the induction hypothesis gives
$$
a_k \leq a_n+a_m \leq C(n+m) = Ck,
$$
which completes the proof.
\end{proof}

\subsection{Convergence in Hilbert spaces}

The next lemma formalizes the computation in Subsection \ref{subsec:defect}. It shows that the distances $d_b$ control the distances between the normalized vectors of an $\A$-admissible sequence in a Hilbert space, whose norms are comparable to the indices.

\begin{lemma}\label{lem:hilbert-path}
Let $\A \subseteq \Delta$, and let $H$ be a real or complex Hilbert space. Further, let $(z_n)_{n \in \mathbb{N}} \subseteq H$ be an $\A$-admissible sequence such that
$$
n \leq \|z_n\| \leq 2n, \quad n \in \N,
$$
an let $b_n = \|z_n\|-n$ and $u_n = z_n/\|z_n\|$ for $n \in \N$. Then, $b \in \Pro(\A)$ and
$$
\|u_r-u_s\| \leq 2d_b(r,s), \quad r,s \in \N.
$$
\end{lemma}

\begin{proof}
Set $a_n=\|z_n\|$. Using admissibility it follows that $D_b=D_a\geq0$ on $\A$. Since $0\leq b_n\leq n$, we also have that $b\in\Pro(\A)$. The computation in Subsection~\ref{subsec:defect} shows that $\|u_n-u_m\|\leq2w_b(n,m)$ on every edge. Summing along a path from $r$ to $s$ and taking the infimum over all such paths gives
$\|u_r-u_s\|\leq2d_b(r,s)$.
\end{proof}

To apply Lemma \ref{lem:hilbert-path} to an arbitrary $\A$-admissible sequence $(v_n)_{n \in \mathbb{N}}$, we adjoin the scalar coordinate $n$ to $v_n$. By Lemma \ref{lem:linear-bound}, the norms of the resulting vectors are comparable to $n$. This yields the main result of this section, which is the content of the following statement.

\begin{proposition}\label{prop:hilbert-sufficiency}
Let $\A \subseteq \Delta$ satisfy condition \eqref{eq:defect-condition} and let $H$ be a real or complex Hilbert space, equipped with the norm $\|x\|=\sqrt{\langle x,x\rangle_H}$. Then $\|\cdot\|$ has the $\A$-Fekete property.
\end{proposition}

\begin{proof}
Let $(v_n)_{n \in \mathbb{N}} \subseteq H$ be an $\A$-admissible sequence. Lemma \ref{lem:linear-bound} shows that there exists $C>0$ such that $\|v_n\| \leq Cn$, for all $n \in \N$. Let $\K$ be the scalar field of $H$ and let $\K \oplus_2 H$ denote the Hilbert direct sum, equipped with the norm
\begin{equation*}
    \|(a,x)\|=\bigl(|a|^2+\|x\|^2\bigr)^{1/2}.
\end{equation*}
Define
$$
z_n = (n,C^{-1}v_n) \in \K \oplus_2 H, \quad n \in \N.
$$
Then, $n \leq \|z_n\| \leq \sqrt{2}\,n$, and for $(n,m) \in \A$ we have
\begin{align*}
\|z_{n+m}\|^2 &= (n+m)^2+C^{-2}\|v_{n+m}\|^2\\
  &\leq (n+m)^2+C^{-2}\|v_n+v_m\|^2\\
  &= \|z_n+z_m\|^2.
\end{align*}
Hence, $(z_n)_{n \in \mathbb{N}}$ is $\A$-admissible.

Let $u_n = z_n/\|z_n\|$. By Lemma \ref{lem:hilbert-path}, there exists $b \in \Pro(\A)$ such that $\|u_r-u_s\| \leq 2d_b(r,s)$ for all $r,s \in \N$. Since $\A$ satisfies condition \eqref{eq:defect-condition}, it follows that $(u_n)_{n \in \mathbb{N}}$ is a Cauchy sequence. Denote by $(t,y) \in \K \oplus_2 H$ the limit of $(u_n)_{n \in \mathbb{N}}$. The first coordinate of $u_n$ is $n/\|z_n\| \in [1/\sqrt{2},1]$, and therefore $t \geq 1/\sqrt{2}$. Consequently,
$$
\frac{v_n}{n} = C\,\frac{C^{-1}v_n/\|z_n\|}{n/\|z_n\|} \longrightarrow \frac{C}{t}\,y, \quad n \to \infty,
$$
which proves the claim with $v = Cy/t$.
\end{proof}

Applying Proposition \ref{prop:hilbert-sufficiency} to the real line yields a scalar version of Fekete's lemma under the defect condition.

\begin{corollary}\label{cor:scalar-convergence}
Let $\A \subseteq \Delta$ satisfy condition \eqref{eq:defect-condition}. Then the following hold.
\begin{enumerate}[\rm (i)]
\item If $a=(a_n)_{n \in \mathbb{N}}$ is an $\A$-subadditive sequence, then the limit $\lim_{n \to \infty} a_n/n$ exists in $[0,\infty)$.
\item If $X$ is a normed space and $(v_n)_{n \in \mathbb{N}} \subseteq X$ is $\A$-admissible, then the limit $\lim_{n \to \infty} \|v_n\|/n$ exists in $[0,\infty)$.
\end{enumerate}
\end{corollary}

\begin{proof}
Since $a$ is nonnegative, we have $|a_{n+m}| = a_{n+m} \leq a_n+a_m = |a_n+a_m|$ for all $(n,m) \in \A$. Thus, $(a_n)_{n \in \mathbb{N}}$ is an $\A$-admissible sequence in $\R$. By Proposition \ref{prop:hilbert-sufficiency}, the limit of $(a_n/n)_{n \in \mathbb{N}}$ exists, and it lies in $[0,\infty)$. This proves the first claim. For the second claim, observe that for $(n,m) \in \A$ we have $\|v_{n+m}\| \leq \|v_n+v_m\| \leq \|v_n\|+\|v_m\|$. Hence, $(\|v_n\|)_{n \in \mathbb{N}}$ is $\A$-subadditive, and the claim follows from (i).
\end{proof}

\section{Necessity of the defect condition}\label{sec:plane}

In this section, we show that if $\A$ does not satisfy condition \eqref{eq:defect-condition}, then every normed space $X$ with $\dim_\R X \geq 2$ contains an $\A$-admissible sequence $(v_n)_{n \in \mathbb{N}}$ such that $(v_n/n)_{n \in \mathbb{N}}$ does not converge (Corollary \ref{cor:plane-obstruction}). It is essential that this holds for every norm. Indeed, a counterexample for the Euclidean norm on $\R^2$ alone would not exclude the existence of another norm on $\R^2$ with the $\A$-Fekete property.

\subsection{A function on the graph}

Recall that a \emph{pseudometric} on a set $M$ is a map $d : M \times M \to [0,\infty)$ such that
$$
d(x,x)=0, \quad d(x,y)=d(y,x), \quad d(x,z) \leq d(x,y)+d(y,z),
$$
for $x, y, z \in M$.
The pair $(M,d)$ is called a \emph{pseudometric space}. A sequence $(x_n)_{n \in \mathbb{N}} \subseteq M$ is called a \emph{Cauchy sequence} if for every $\varepsilon>0$ there exists $N \in \N$ such that $d(x_n,x_m)<\varepsilon$ for all $n,m \geq N$. The next lemma shows that if a sequence is not a Cauchy sequence, then there exists a bounded $1$-Lipschitz function whose values along the sequence do not converge.

\begin{lemma}\label{lem:hilbert-metric-test}
Let $(M,d)$ be a pseudometric space such that $d(x,y) \leq 1$ for all $x,y \in M$, and let $(x_n)_{n \in \mathbb{N}} \subseteq M$ be a sequence which is not a Cauchy sequence. Then, there exists a function $f : M \to [0,1]$ such that
$$
|f(x)-f(y)| \leq d(x,y), \quad x,y \in M,
$$
for which the sequence $(f(x_n))_{n \in \mathbb{N}}$ does not converge.
\end{lemma}

\begin{proof}

Since $(x_n)_{n \in \mathbb{N}}$ is not Cauchy, there exist $\varepsilon>0$ and subsequences
$(u_j)_{j \in \mathbb{N}}$ and $(v_j)_{j \in \mathbb{N}}$ of $(x_n)_{n \in \mathbb{N}}$ such that $d(u_j,v_j)\geq\varepsilon$
for all $j\in\N$. Applying Efremovi\v{c}'s lemma
\cite[Lemma~2.2]{kunduJain2006} in the metric quotient of $M$, and
passing to subsequences $(u_{j_k})_{k\in\N}$ and $(v_{j_k})_{k\in\N}$
along a common sequence of indices $j_1<j_2<\cdots$, which we again
denote by $(u_j)_{j \in \mathbb{N}}$ and $(v_j)_{j \in \mathbb{N}}$, we obtain
$$
d(u_i,v_j)\geq\varepsilon/4, \quad i,j\in\N.
$$
We now define $f(x)=\inf_j d(x,u_j)$. Then, $f$ is $1$-Lipschitz and takes
values in $[0,1]$. Since $f(u_j)=0$ and $f(v_j)\geq\varepsilon/4$
for every $j\in\N$, the sequence $(f(x_n))_{n \in \mathbb{N}}$ does not converge.
\end{proof}

\subsection{A quadratic estimate on the unit sphere}

Let $Y$ be a real normed space with $\dim Y=2$. We denote its unit sphere by $S_Y = \{y \in Y : \|y\|=1\}$. For a linear functional $F : Y \to \R$ we write $\|F\| = \sup_{y \in S_Y} |F(y)|$. If $Y$ is the Euclidean plane, then the functions $\gamma(q) = (\cos q,\sin q)$ and $F_q = \langle \cdot,\gamma(q) \rangle$ satisfy $F_q(\gamma(q))=1$, as well as
$$
0 \leq 1-F_q(\gamma(q')) = 1-\cos(q-q') \leq \frac{(q-q')^2}{2}, \quad q,q' \in \R.
$$
The following lemma shows that an estimate of this type holds along a short arc of $S_Y$ for every norm on $Y$.

\begin{lemma}\label{lem:hilbert-supporting-arc}
Let $Y$ be a real normed space with $\dim Y=2$ and let $\eta>0$. Then, there exist real numbers $q_-<q_+$ with $q_+-q_-<\eta$, a continuous map $\gamma : [q_-,q_+] \to S_Y$ with $\gamma(q_-) \neq \gamma(q_+)$, and a family $\{F_q\}_{q \in [q_-,q_+]}$ of linear functionals $F_q : Y \to \R$ with $\|F_q\|=1$, such that
\begin{equation}\label{eq:hilbert-arc-defect}
F_q(\gamma(q))=1, \quad 0 \leq 1-F_q(\gamma(q')) \leq \frac{(q-q')^2}{4}, \quad q,q' \in [q_-,q_+].
\end{equation}
\end{lemma}

\begin{proof}
We identify $Y$ with $\R^2$ and denote by $K = \{z \in \R^2 : \|z\| \leq 1\}$ the closed unit ball of $Y$. Since all norms on $\R^2$ are equivalent, we can choose a Euclidean inner product $\langle \cdot,\cdot \rangle$, with associated norm $|\cdot|$, such that the closed Euclidean unit disk is contained in $K$. We fix orthonormal coordinates with respect to this inner product and set
$$
R = 1+\max\{x : (x,y) \in K\}, \quad z(q) = (R,q), \quad q \in \R.
$$

\smallskip
\textbf{Step 1. Construction of $\gamma$ and $F_q$.}
Let $P_K$ denote the Euclidean projection onto $K$ and set
$\gamma(q)=P_Kz(q)$ and $\nu(q)=z(q)-\gamma(q)$. By
\cite[Theorem~5.2]{brezis2011} and \cite[Proposition~5.3]{brezis2011}, $P_K$ is
$1$-Lipschitz and
\begin{equation}\label{eq:projection-support}
\langle \nu(q),z-\gamma(q)\rangle\leq 0, \quad z\in K.
\end{equation}
In particular, $\gamma$ is continuous. Since $K$ contains the
Euclidean unit disk and the first coordinate of $\nu(q)$ is at least
one, we have
$$
h(q):=\langle\nu(q),\gamma(q)\rangle
     =\sup_{z\in K}\langle\nu(q),z\rangle\geq|\nu(q)|\geq1.
$$
Defining $F_q(z)=\langle\nu(q),z\rangle/h(q)$ and using the symmetry of $K$,
we have $\|F_q\|=1$. Since $F_q(\gamma(q))=1$ and $\gamma(q)\in K$,
it follows that $\gamma(q)\in S_Y$.

\smallskip
\textbf{Step 2. The quadratic estimate.}
Let $u=\gamma(q)-\gamma(q')$ and $v=\nu(q)-\nu(q')$. Applying
\eqref{eq:projection-support} at $q$ with $z=\gamma(q')$ and at $q'$
with $z=\gamma(q)$ gives
$0\leq\langle\nu(q),u\rangle\leq\langle v,u\rangle$, and therefore
\begin{align*}
0 &\leq 1-F_q(\gamma(q'))\\
  &=\frac{\langle\nu(q),u\rangle}{h(q)}\\
  &\leq \langle v,u\rangle\\
  &=\frac{|u+v|^2-|u-v|^2}{4}\\
  &\leq\frac{(q-q')^2}{4},
\end{align*}
where we used that $u+v=z(q)-z(q')=(0,q-q')$.

\smallskip
\textbf{Step 3. Choice of $q_-$ and $q_+$.} We proceed by showing that $\gamma$ is not constant and argue by contradiction. Suppose by contradiction that $\gamma(q) = y_0$ for all $q \in \R$. Then \eqref{eq:projection-support} gives
$$
\langle z(q)-y_0,z-y_0 \rangle \leq 0, \quad z \in K, \ q \in \R.
$$
Dividing by $|q|$ and letting $q \to \pm\infty$ shows that every point of $K$ has the same second coordinate as $y_0$. This contradicts the inclusion of the Euclidean unit disk in $K$. Hence, there exist $t_-<t_+$ with $\gamma(t_-) \neq \gamma(t_+)$.

We subdivide $[t_-,t_+]$ into finitely many intervals of length less than $\eta$. Since $\gamma(t_-) \neq \gamma(t_+)$, the map $\gamma$ takes distinct values at the endpoints $q_-<q_+$ of one of these intervals. Restricting $\gamma$ and $(F_q)$ to $[q_-,q_+]$ proves the lemma.
\end{proof}

\subsection{Construction of the counterexample}

The next proposition is the main result of this section. It turns a sequence $b \in \Pro(\A)$ which witnesses the failure of condition \eqref{eq:defect-condition} into an $\A$-admissible sequence whose directions do not converge.

\begin{proposition}\label{prop:plane-obstruction}
Let $\A \subseteq \Delta$, and let $b \in \Pro(\A)$ be such that
\begin{equation}\label{eq:failure}
\lim_{N \to \infty} \sup_{r,s \geq N} d_b(r,s)>0.
\end{equation}
Further, let $X$ be a normed space and let $Y \subseteq X$ be a two-dimensional real subspace. Then, there exists an $\A$-admissible sequence $(v_n)_{n \in \mathbb{N}} \subseteq Y$ satisfying
$$
\|v_n\| = n+b_n, \quad n \in \N,
$$
and such that $(v_n/\|v_n\|)_{n \in \mathbb{N}}$ has two distinct accumulation points. Moreover, the sequence $(v_n/n)_{n \in \mathbb{N}}$ does not converge weakly in $X$.
\end{proposition}

\begin{proof}
\textbf{Step 1. A function on $\N$.} For $r,s \in \N$, define the map $\bar d(r,s) = \min\{1,d_b(r,s)\}$. By \eqref{eq:pseudometric} and the inequality $\min\{1,x+y\} \leq \min\{1,x\}+\min\{1,y\}$, the map $\bar d$ is a pseudometric on $\N$ with values in $[0,1]$.

Since the supremum in \eqref{eq:failure} is nonincreasing in $N$, there exists $\varepsilon \in (0,1)$ such that $\sup_{r,s \geq N} d_b(r,s)>\varepsilon$ for every $N \in \N$. Hence, for every $N \in \N$ there exist $r,s \geq N$ with $\bar d(r,s)>\varepsilon$. In other words, the sequence $(n)_{n \in \N}$ is not a Cauchy sequence in $(\N,\bar d)$. By Lemma \ref{lem:hilbert-metric-test}, there exists a function $f : \N \to [0,1]$ such that $(f(n))_{n \in \mathbb{N}}$ does not converge and
$$
|f(n)-f(m)| \leq \bar d(n,m), \quad n,m \in \N.
$$
Let $\alpha = \liminf_{n \to \infty} f(n)$ and $\beta = \limsup_{n \to \infty} f(n)$. Since $(f(n))_{n \in \mathbb{N}}$ does not converge, we have $\rho = \beta-\alpha>0$. Further, define $g_n$ via
$$
g_n = \min\Big\{1,\max\Big\{0,\frac{f(n)-\alpha}{\rho}\Big\}\Big\}, \quad n \in \N.
$$
Since the map $x \mapsto \min\{1,\max\{0,x\}\}$ is $1$-Lipschitz, we have
\begin{equation}\label{eq:g-lipschitz}
0 \leq g_n \leq 1, \quad |g_n-g_m| \leq \rho^{-1}\,\bar d(n,m), \quad n,m \in \N.
\end{equation}
Moreover, $g_n \to 0$ along every subsequence for which $f(n) \to \alpha$, and $g_n \to 1$ along every subsequence for which $f(n) \to \beta$.

\smallskip
\textbf{Step 2. Definition of $(v_n)_{n \in \mathbb{N}}$.} We apply Lemma \ref{lem:hilbert-supporting-arc} to $Y$ and $\eta = \sqrt{2}\,\rho$, and we let $q_-,q_+$, $\gamma$ and $(F_q)_{q \in [q_-,q_+]}$ be given as in the lemma. Set $L = q_+-q_-<\sqrt{2}\,\rho$ and define
$$
q_n = q_-+Lg_n, \quad a_n = n+b_n, \quad v_n = a_n\gamma(q_n), \quad n \in \N.
$$
Then, $q_n \in [q_-,q_+]$, $n \leq a_n \leq 2n$ and $\|v_n\| = a_n = n+b_n$. For $n,m \in \N$, the estimates \eqref{eq:hilbert-arc-defect} and \eqref{eq:g-lipschitz} yield
\begin{equation}\label{eq:hilbert-plane-directed}
0 \leq 1-F_{q_m}(\gamma(q_n)) \leq \frac{L^2(g_n-g_m)^2}{4} \leq \frac{L^2}{4\rho^2}\,\bar d(n,m)^2 \leq \frac{1}{2}\,\bar d(n,m)^2.
\end{equation}

\smallskip
\textbf{Step 3. Admissibility.} Let $(n,m) \in \A$ and $D = D_b(n,m)$. Since $a_n = n+b_n$, we have $D = D_a(n,m)$, and \eqref{eq:one-edge} gives $\bar d(n,m)^2 \leq d_b(n,m)^2 \leq D/n$. Using $\|F_{q_m}\|=1$, $F_{q_m}(\gamma(q_m))=1$ as well as the estimate \eqref{eq:hilbert-plane-directed} and $a_n \leq 2n$, we obtain
$$
\begin{aligned}
\|v_n+v_m\| &\geq F_{q_m}(v_n+v_m) = a_nF_{q_m}(\gamma(q_n))+a_m \\
  &\geq a_n+a_m-\frac{a_n}{2}\,\bar d(n,m)^2 \\
&\geq a_n+a_m-D = a_{n+m} = \|v_{n+m}\|.
\end{aligned}
$$
This shows that $(v_n)_{n \in \mathbb{N}}$ is $\A$-admissible.

\smallskip
\textbf{Step 4. Non-convergence.} By Step 1 and the continuity of $\gamma$, the sequence $v_n/\|v_n\| = \gamma(q_n)$ has the distinct accumulation points $\gamma(q_-)$ and $\gamma(q_+)$. Since $a_n/n \in [1,2]$, we may pass to further subsequences along which $a_n/n$ converges. Thus, $(v_n/n)_{n \in \mathbb{N}}$ has accumulation points $r\gamma(q_-)$ and $s\gamma(q_+)$ for some $r,s \in [1,2]$. Moreover, we have $\gamma(q_-) \neq \gamma(q_+)$.

If $(v_n/n)_{n \in \mathbb{N}}$ converged weakly to some $v \in X$, then $v$ would coincide with both accumulation points $r\gamma(q_-)$ and $s\gamma(q_+)$, since norm convergence implies weak convergence and weak limits are unique. Consequently, $(v_n/n)_{n \in \mathbb{N}}$ does not converge weakly in $X$.
\end{proof}

Since the failure of condition \eqref{eq:defect-condition} provides a sequence $b \in \Pro(\A)$ satisfying \eqref{eq:failure}, Proposition \ref{prop:plane-obstruction} has the following consequence.

\begin{corollary}\label{cor:plane-obstruction}
Let $\A \subseteq \Delta$ be such that condition \eqref{eq:defect-condition} fails and let $X$ be a normed space with $\dim_\R X \geq 2$. Then, there exists an $\A$-admissible sequence $(v_n)_{n \in \mathbb{N}} \subseteq X$ satisfying
$$
n \leq \|v_n\| \leq 2n, \quad n \in \N,
$$
and such that $(v_n/n)_{n \in \mathbb{N}}$ does not converge weakly. In particular, the norm of $X$ does not have the $\A$-Fekete property.
\end{corollary}

\begin{proof}
Since $\A$ does not satisfy condition \eqref{eq:defect-condition}, there exists $b \in \Pro(\A)$ so that \eqref{eq:failure} holds. Let $Y \subseteq X$ be a real subspace of dimension $2$. Proposition \ref{prop:plane-obstruction} yields an $\A$-admissible sequence $(v_n)_{n \in \mathbb{N}} \subseteq Y$ with $\|v_n\| = n+b_n \in [n,2n]$ such that $(v_n/n)_{n \in \mathbb{N}}$ does not converge weakly. Since norm convergence implies weak convergence, the claim follows.
\end{proof}

We are now ready to prove Theorem \ref{thm:test-spaces}.

\begin{proof}[Proof of Theorem \ref{thm:test-spaces}]
The implication (i)$\Rightarrow$(iii) is Proposition \ref{prop:hilbert-sufficiency}. The implication (iii)$\Rightarrow$(ii) holds since the Euclidean plane is a real Hilbert space. Finally, the implication (ii)$\Rightarrow$(i) follows from Corollary \ref{cor:plane-obstruction}, applied to the Euclidean plane.
\end{proof}

\section{Reflexivity and an equivalent Fekete norm}\label{sec:renorming}

The goal of this section is to prove Theorem \ref{thm:reflexivity-characterization}. To that end, let $X$ be a separable Banach space with $\dim X \geq 2$ and let $\A \subseteq \Delta$. If $X$ admits an equivalent norm $q$ with the $\A$-Fekete property, then Corollary \ref{cor:plane-obstruction}, applied to the normed space $(X,q)$, shows that $\A$ satisfies condition \eqref{eq:defect-condition}. In Subsection \ref{subsec:james}, we show that $X$ is necessarily reflexive. Conversely, in Subsections \ref{subsec:rows}--\ref{subsec:proof-thm13} we show that every separable reflexive Banach space admits an equivalent norm $\mathfrak{q}$ which has the $\A$-Fekete property for every $\A \subseteq \Delta$ satisfying condition \eqref{eq:defect-condition} (Theorem \ref{thm:sufficiency}).

\subsection{Necessity of reflexivity}\label{subsec:james}

Let $X$ be a Banach space over $\K$. We denote by $X_\R$ the real Banach space obtained from $X$ by restricting the scalar multiplication to $\R$. Recall that, by Kakutani's theorem, a Banach space is reflexive if and only if its closed unit ball is weakly compact. If $X$ is complex, then every continuous real-linear functional $f : X \to \R$ is the real part of the continuous complex-linear functional $x \mapsto f(x)-\mathrm if(\mathrm ix)$. Consequently, $X$ and $X_\R$ carry the same weak topology. Moreover, the weak topology does not change if the norm of $X$ is replaced by an equivalent norm $q$. The closed unit balls of $q$ and of the original norm are weakly closed, and each of them is contained in a multiple of the other. Hence, if $X$ is not reflexive, then $X_\R$, equipped with $q$, is not reflexive either.

The next lemma shows that in a nonreflexive space the $\A$-Fekete property fails for every equivalent norm, even if norm convergence is replaced by weak convergence. It does not require the Banach space $X$ to be separable.

\begin{lemma}\label{lem:james-obstruction}
Let $(X,\|\cdot\|)$ be a Banach space which is not reflexive and let $q$ be a norm on $X$ which is equivalent to $\|\cdot\|$. Then, there exists a sequence $(v_n)_{n \in \mathbb{N}} \subseteq X$ satisfying
$$
q(v_{n+m})<q(v_n+v_m), \quad n,m \in \N,
$$
and such that $(v_n/n)_{n \in \mathbb{N}}$ does not converge weakly. In particular, $q$ does not have the $\A$-Fekete property for any $\A \subseteq \Delta$.
\end{lemma}

\begin{proof}
Since the Banach space $(X_\R,q)$ is not reflexive, James's theorem \cite{james1972reflexivity} gives a continuous linear functional $f : X \to \R$ of norm one that does not attain its supremum on the closed unit ball of $q$, i.e.,
$$
\sup_{q(x) \leq 1} f(x) = 1, \quad f(x)<1 \ \text{ whenever } q(x) \leq 1.
$$
In particular, $f(x) \leq q(x)$ for all $x \in X$.

Let $c_n = 1+2^{-n}$. Since $c_{n+1}/c_n<1$, we can choose $x_n \in X$ with $q(x_n)=1$ and $f(x_n)>c_{n+1}/c_n$. Further, define $v_n = nc_nx_n$. Since $(c_n)_{n \in \mathbb{N}}$ is decreasing, for all $n,m \in \N$ we have
$$
\begin{aligned}
q(v_n+v_m) &\geq f(v_n+v_m) = nc_nf(x_n)+mc_mf(x_m) \\
&> nc_{n+1}+mc_{m+1} \geq (n+m)c_{n+m} = q(v_{n+m}).
\end{aligned}
$$
In particular, $(v_n)_{n \in \mathbb{N}}$ is $\A$-admissible with respect to $q$ for every $\A \subseteq \Delta$.

It remains to show that $(v_n/n)_{n \in \mathbb{N}}$ does not converge weakly. To do so, we suppose by contradiction that $v_n/n = c_nx_n$ converges weakly to some $v \in X$. Since $q$ is weakly lower semicontinuous, we have $q(v) \leq \lim_{n \to \infty} c_n = 1$. On the other hand, $c_{n+1}<f(c_nx_n) \leq c_n$, and the weak continuity of $f$ gives $f(v)=1$. This contradicts the fact that $f(x)<1$ whenever $q(x) \leq 1$.
\end{proof}

\subsection{Rows with infinitely many small defects}\label{subsec:rows}

Let $\A \subseteq \Delta$. For $n \in \N$, we call the set $\{m \in \N : (n,m) \in \A\}$ the \emph{$n$-th row} of $\A$. The defect condition allows paths to pass through small indices, and it imposes no direct bound on the costs of individual edges. Nevertheless, it implies that rows with arbitrarily large index contain infinitely many edges with small defect (see Lemma \ref{lem:small-defect-row}). This is a consequence of the next lemma, which is also used in Section \ref{sec:consequences}. It shows that if every row with large index contains only finitely many edges with small defect, then the defect condition fails for an explicit sequence in $\Pro(\A)$.

\begin{lemma}\label{lem:block-construction}
Let $\A \subseteq \Delta$ and let $a=(a_n)_{n \in \mathbb{N}}$ be an $\A$-subadditive sequence with $a_n \leq n$ for all $n \in \N$. Further, let $\varepsilon>0$ and $K \in \N$, and assume that for every $n \geq K$ the set
$$
R_n = \{m \in \N : (n,m) \in \A, \ D_a(n,m)<\varepsilon n\},
$$
is finite. Then, there exists a nondecreasing function $\kappa : \N \to \N$ with $\kappa(n) \to \infty$ as $n \to \infty$ such that the sequence $b=(b_n)_{n \in \mathbb{N}}$ defined by
$$
b_n = \frac{1}{3}\Big(a_n+\frac{n}{\kappa(n)}+1\Big), \quad n \in \N,
$$
belongs to $\Pro(\A)$ and satisfies \eqref{eq:failure}.
\end{lemma}

\begin{proof}We divide the proof into four main steps.

\textbf{Step 1: A decomposition into blocks.} For $n \geq K$, let $g(n) = \max(\{n\} \cup R_n)$, which is well defined because $R_n$ is finite. Set $M_1 = K$ and define the numbers $M_{j+1}$ recursively via
$$
M_{j+1} = 1+\max\big(\{M_j\} \cup \{g(n) : K \leq n<M_j\}\big), \quad j \in \N.
$$
Then, $(M_j)_{j \in \mathbb{N}}$ is strictly increasing, and $g(n)<M_{j+1}$ whenever $K \leq n<M_j$. Next, define
\begin{equation}
\kappa(n) =
\begin{cases}
j, & \text{if } M_j \leq n<M_{j+1}, \\
1, & \text{if } n<K.
\end{cases}
\end{equation}
This way, $\kappa$ is nondecreasing, tends to infinity and satisfies $\kappa(M_j)=j$. We call the set $\kappa^{-1}(\{j\})$ the $j$-th block.

Now, let $(n,m) \in \A$ with $n \geq K$ and $D_a(n,m)<\varepsilon n$. Writing $j = \kappa(n)$, we have $n<M_{j+1}$, and hence $n \leq m \leq g(n)<M_{j+2}$. Since $\kappa$ is nondecreasing, this shows that for all $(n,m) \in \A$ we have the implication
\begin{equation}\label{eq:consecutive-blocks}
n \geq K \ \text{ and } \ D_a(n,m)<\varepsilon n \quad \Longrightarrow \quad \kappa(m) \in \{\kappa(n),\kappa(n)+1\}.
\end{equation}
In other words, an edge with small defect whose smaller endpoint is at least $K$ joins vertices in the same block or in consecutive blocks.

\smallskip
\textbf{Step 2: The sequence $b$.} Set $c_n = n/\kappa(n)$. Since $\kappa$ is nondecreasing, we have
$$
c_{n+m} = \frac{n+m}{\kappa(n+m)} \leq \frac{n}{\kappa(n)}+\frac{m}{\kappa(m)} = c_n+c_m, \quad n,m \in \N.
$$
Thus, $D_c(n,m) \geq 0$ for all $n,m \in \N$. Since $b_n = (a_n+c_n+1)/3$ and $a$ is $\A$-subadditive, it follows that
\begin{equation}\label{eq:block-defect}
D_b(n,m) = \frac{D_a(n,m)+D_c(n,m)+1}{3} \geq \frac{1}{3}, \quad (n,m) \in \A.
\end{equation}
Moreover, $0 \leq b_n \leq (2n+1)/3 \leq n$, and hence $b \in \Pro(\A)$.

\smallskip
\textbf{Step 3: Lower bounds for edge costs.} Recall from \eqref{eq:one-edge} that an edge $(n,m) \in \A$ has cost $w_b(n,m) = (D_b(n,m)/n)^{1/2}$. By \eqref{eq:block-defect}, its cost satisfies
\begin{equation}
w_b(n,m) \geq
\begin{cases}
\sqrt{\varepsilon/3}, & \text{if } D_a(n,m) \geq \varepsilon n, \\
1/\sqrt{3K}, & \text{if } n<K.
\end{cases}
\end{equation}
Now let $(n,m) \in \A$ with $\kappa(n)=j$ and $\kappa(m)=j+1$. Then $\kappa(n+m) \geq j+1$, and hence
$$
D_c(n,m) \geq \frac{n}{j}+\frac{m}{j+1}-\frac{n+m}{j+1} = \frac{n}{j(j+1)}.
$$
Consequently, by \eqref{eq:block-defect}, every edge joining the $j$-th and the $(j+1)$-th block has cost at least $1/\sqrt{3j(j+1)}$.

\smallskip
\textbf{Step 4: Large distances.} Choose $0<\eta<\min\{\sqrt{\varepsilon/3},1/\sqrt{3K}\}$. Since the series $\sum_i (3i(i+1))^{-1/2}$ diverges, it follows that for every $j \in \N$ there exists $k>j$ such that
$$
\sum_{i=j}^{k-1} \frac{1}{\sqrt{3i(i+1)}} \geq \eta.
$$
We claim that $d_b(M_j,M_k) \geq \eta$. To see this, consider a path from $M_j$ to $M_k$. If one of its edges has cost at least $\eta$, then the total cost of the path is at least $\eta$ as well. Otherwise, Step 3 shows that every edge $(n,m) \in \A$ of the path satisfies $n \geq K$ and $D_a(n,m)<\varepsilon n$. By \eqref{eq:consecutive-blocks}, every edge of the path then joins vertices in the same block or in consecutive blocks. Since $\kappa(M_j)=j$ and $\kappa(M_k)=k$, the path contains for every $i \in \{j,\dots,k-1\}$ an edge joining the $i$-th and the $(i+1)$-th block. These edges are distinct, and by Step 3 their costs sum to at least $\eta$. This proves the claim. Since $M_j \to \infty$ as $j \to \infty$, it follows that $\sup_{r,s \geq N} d_b(r,s) \geq \eta$ for every $N \in \N$. Consequently, $b$ satisfies \eqref{eq:failure}.
\end{proof}

We now deduce the property of the rows of $\A$ mentioned at the beginning of this subsection: if $\A$ satisfies condition \eqref{eq:defect-condition}, then rows of $\A$ with arbitrarily large index contain infinitely many edges with small defect.

\begin{lemma}\label{lem:small-defect-row}
Let $\A \subseteq \Delta$ satisfy condition \eqref{eq:defect-condition} and let $a=(a_n)_{n \in \mathbb{N}}$ be an $\A$-subadditive sequence for which there exists $C>0$, such that $a_n \leq Cn$ for all $n \in \N$. Then, for every $\varepsilon>0$ and every $K \in \N$ there exists $n \geq K$ such that the set
\begin{equation}\label{eq:infinite-small-defect-row}
\{m \in \N : (n,m) \in \A, \ D_a(n,m)<\varepsilon n\},
\end{equation}
is infinite.
\end{lemma}

\begin{proof}
Suppose by contradiction that there exist $\varepsilon>0$ and $K \in \N$ such that the set in equation \eqref{eq:infinite-small-defect-row} is finite for every $n \geq K$. The sequence $a/C$ is $\A$-subadditive and satisfies $a_n/C \leq n$ for all $n \in \N$. Moreover, $D_{a/C}(n,m)<(\varepsilon/C)n$ if and only if $D_a(n,m)<\varepsilon n$. Hence, Lemma \ref{lem:block-construction}, applied to $a/C$, $\varepsilon/C$ and $K$, yields a sequence in $\Pro(\A)$ for which condition \eqref{eq:defect-condition} fails. This contradicts the assumption on $\A$.
\end{proof}

\subsection{\texorpdfstring{$2$}{2}-rotund norms}\label{subsec:2rotund}

The construction of the norm $\mathfrak{q}$ starts from a norm with the following convexity property.

\begin{definition}
Let $X$ be a vector space over $\K$. A norm $p$ on $X$ is called \emph{$2$-rotund} if, every sequence $(x_n)_{n \in \mathbb{N}} \subseteq X$ with $p(x_n)=1$ for all $n \in \N$ and
$$
\lim_{n,m \to \infty} p(x_n+x_m) = 2,
$$
converges in $(X,p)$. 
\end{definition}

The next lemma, which is a consequence of a renorming theorem of Odell and Schlumprecht \cite{odell1998asymptotic}, shows that such a norm exists on every separable reflexive space.

\begin{lemma}\label{lem:2rotund}
Let $X$ be a separable reflexive Banach space over $\K$. Then, there exists a $2$-rotund norm on $X$ which is equivalent to the norm of $X$.
\end{lemma}

\begin{proof} We divide the proof into two main steps.

\textbf{Step 1: The real case.} Suppose that $\K=\R$. Since $X$ is
separable and reflexive, the renorming theorem of Odell and
Schlumprecht \cite[Corollary~3.1(a)]{odell1998asymptotic} gives an
equivalent $2$-rotund norm $p_0$ on $X$. For the formulation with
the joint limit $\lim_{n,m\to\infty}p_0(x_n+x_m)$, see the paragraph
following Definition~1.1 in \cite{dilworth2024}, where the slightly different definitions occurring in \cite{odell1998asymptotic,dilworth2024} are shown equivalent.

\smallskip
\textbf{Step 2: The complex case.} Suppose now that $\K=\C$. We apply Step 1 to the separable reflexive real Banach space $X_\R$ and let $p_0$ be the resulting real norm. Choose $0<c \leq C<\infty$ with $c\|x\| \leq p_0(x) \leq C\|x\|$, and define
$$
p(x) = \Big(\int_\T p_0(\lambda x)^2\,d\mu(\lambda)\Big)^{1/2}, \quad x \in X,
$$
where $\mu$ is the normalized arc length measure on $\T$. Minkowski's inequality and the rotation invariance of $\mu$ show that $p$ is a complex norm. Moreover, $c\|x\| \leq p(x) \leq C\|x\|$ for all $x \in X$.

Now, suppose that $p(x_n)=1$ for all $n\in\N$ and that $p(x_n+x_m)\to2$
as $n,m\to\infty$. For $\lambda\in\T$, let
$$
h_{n,m}(\lambda)
=2p_0(\lambda x_n)^2+2p_0(\lambda x_m)^2
 -p_0(\lambda(x_n+x_m))^2.
$$
By the triangle inequality, we have $h_{n,m}\geq0$. Moreover, $(x_n)_{n \in \mathbb{N}}$ is
bounded, so the functions $h_{n,m}$ have a common Lipschitz constant.

Since
$$
\int_\T h_{n,m}\,d\mu=4-p(x_n+x_m)^2\to0,
$$
it follows that $h_{n,m}(1)\to0$. Indeed, a positive lower bound at
$1$ along a sequence of pairs would give a positive lower bound on
a fixed arc about $1$. Applying the characterization of $2$-rotundity
due to Dilworth, Kutzarova and Motakis \cite[Section~2, (2.1)]{dilworth2024} to the real norm $p_0$, we obtain that
$(x_n)_{n \in \mathbb{N}}$ converges in norm. Hence, $p$ is $2$-rotund.
\end{proof}

The next lemma is the consequence of $2$-rotundity which we use in the proof of Theorem \ref{thm:sufficiency}.

\begin{lemma}\label{lem:translated-limits}
Let $X$ be a vector space over $\K$ and let $p$ be a $2$-rotund norm on $X$. Further, let $(y_n)_{n \in \mathbb{N}} \subseteq X$ satisfy $p(y_n)=1$ for all $n \in \N$, and assume that the limit
$$
\tau(x) = \lim_{m \to \infty} p(x+y_m)
$$
exists for every $x \in X$. If $(y_n)_{n \in \mathbb{N}}$ does not converge in $(X,p)$, then
$$
\limsup_{n \to \infty} \tau(y_n)<2.
$$
\end{lemma}

\begin{proof}
Suppose that $\limsup_n\tau(y_n)=2$. Since each limit
$\lim_m p(y_n+y_m)=\tau(y_n)$ exists, we may choose a subsequence
$(z_j)_{j \in \mathbb{N}}$ of $(y_n)_{n \in \mathbb{N}}$ such that $\tau(z_j)>2-1/j$, and
$$
|p(z_i+z_j)-\tau(z_i)|<1/j,\quad i<j.
$$
Thus, $2\geq p(z_i+z_j)>2-1/i-1/j$ for $i<j$. Since $p(2z_i)=2$,
it follows from the property of $p$ being $2$-rotund that $z_j\to z$ for some $z\in X$.
Evaluating the limit $\tau(-z)=\lim_m p(y_m-z)$ along $(z_j)_{j \in \mathbb{N}}$
gives $\tau(-z)=0$. Hence $y_m\to z$, which gives a contradiction. Since
$\tau(y_n)\leq2$ for all $n$, the claim follows.
\end{proof}

\subsection{Recording translated norms in a Hilbert space}\label{subsec:translated-norms}

In this subsection, we construct the norm $\mathfrak{q}$ used in the proof of Theorem \ref{thm:sufficiency}. The idea is to record the norms of the translates of a vector in countably many directions as coordinates of a vector in a Hilbert space.

Throughout this subsection, let $X$ be a Banach space over $\K$ and let $p$ be a norm on $X$ which is equivalent to the norm of $X$. Let $U=\{-1,1\}$ if $\K=\R$, and $U=\T$ if $\K=\C$. Further, let $\mu$ be the normalized Haar measure on $U$, i.e., $\mu(\{-1\})=\mu(\{1\})=\frac12$ if $\K=\R$, and $\mu$ is the normalized arc length measure on $\T$ if $\K=\C$. Recall that $\mu$ is invariant under multiplication by elements of $U$, i.e., for every $\omega \in U$ and every continuous function $\phi : U \to \R$ we have
$$
\int_U \phi(\omega\lambda)\,d\mu(\lambda) = \int_U \phi(\lambda)\,d\mu(\lambda).
$$
We denote by $L^2(U)$ the real Hilbert space of square-integrable real-valued functions on $(U,\mu)$.

For $c \in X$, define $\Phi_c : X \to L^2(U)$ and $q_c : X \to [0,\infty)$ by
$$
\Phi_c(x)(\lambda) = p(x+\lambda\,p(x)\,c), \quad q_c(x) = \|\Phi_c(x)\|_{L^2(U)}.
$$
Notice that $\Phi_c(x)$ is a continuous function on $U$. It encodes the norms of the translates of $x$ in the directions $\lambda c$, $\lambda \in U$, scaled by $p(x)$. The averaging over $U$ ensures homogeneity. The next lemma collects the properties of $\Phi_c$ and $q_c$ which we need in the remainder of this section.

\begin{lemma}\label{lem:qc}
Let $X$, $p$, $U$ and $\mu$ be as above and let $c \in X$. Then, $q_c$ is a norm on $X$. Moreover, for all $x,y \in X$ and $t \geq 0$ the following hold:
\begin{enumerate}[\rm (i)]
\item $p(x) \leq q_c(x) \leq (1+p(c))\,p(x)$,
\item $q_c(x+y) \leq \|\Phi_c(x)+\Phi_c(y)\|_{L^2(U)} \leq q_c(x)+q_c(y)$,
\item $\Phi_c(tx) = t\,\Phi_c(x)$.
\end{enumerate}
\end{lemma}

\begin{proof}
For $z\in X$, define $G_z:[0,\infty)\to\R$ via
$$
G_z(s)=\int_U p(z+\lambda sc)^2\,d\mu(\lambda).
$$
By \cite[Lemma~1.10]{maligranda2008}, the sum
$p(z+\lambda sc)^2+p(z-\lambda sc)^2$ is nondecreasing in $s$.
Integrating in $\lambda$ and using the symmetry of $\mu$, we see
that $G_z$ is nondecreasing. Thus, we have
$$
p(x)^2=G_x(0)\leq G_x(p(x))=q_c(x)^2.
$$
The upper bound in (i) follows from the triangle inequality.
For (ii), let $z=x+y$. Then
$$
q_c(z)^2\leq G_z(p(x)+p(y))
\leq\int_U\big(\Phi_c(x)(\lambda)+\Phi_c(y)(\lambda)\big)^2\,d\mu(\lambda),
$$
where the second inequality follows from the triangle inequality for
$p$. Minkowski's inequality gives the remaining estimate in (ii).
Assertion (iii) follows from positive homogeneity of $p$.

Finally, if $\zeta\in\K\setminus\{0\}$, then
$\Phi_c(\zeta x)(\lambda)=|\zeta|\Phi_c(x)(\lambda|\zeta|/\zeta)$.
By the invariance of $\mu$, it follows that $q_c(\zeta x)=|\zeta|q_c(x)$, and together
with (i) and (ii), this shows that $q_c$ is a norm.
\end{proof}

Now let $(c_j)_{j \in \N}$ be a sequence in $X$ and define
$$
w_j = \frac{2^{-j}}{(1+p(c_j))^2}, \quad \mathfrak{q}(x) = \Big(p(x)^2+\sum_{j=1}^\infty w_j\,q_{c_j}(x)^2\Big)^{1/2}, \quad x \in X.
$$
Further, let $H$ be the real Hilbert space
$$
H = \R \oplus_2 \Big(\bigoplus_{j \in \N} L^2(U)\Big)_{\ell^2}.
$$
Its elements are sequences $(t,(\phi_j))$ with $t \in \R$ and $\phi_j \in L^2(U)$ for which the squared norm $\|\cdot\|_H$ satisfies
$$
\|(t,(\phi_j))\|_H^2 = t^2+\sum_{j=1}^\infty \|\phi_j\|_{L^2(U)}^2 < \infty.
$$
Moreover, define the map $F : X \to H$ by
$$
F(x) = \big(p(x),(\sqrt{w_j}\,\Phi_{c_j}(x))_{j \in \N}\big).
$$
The next lemma shows that $\mathfrak{q}$ is a norm that is equivalent to $p$, and that the triangle inequality for $\mathfrak{q}$ is inherited from $H$ via the map $F$.

\begin{lemma}\label{lem:N-F}
Let $X$, $p$, $(c_j)_{j \in \mathbb{N}}$, $\mathfrak{q}$, $H$ and $F$ be as above. Then $F$ is well defined and $\mathfrak{q}$ is a norm on $X$ satisfying $p \leq \mathfrak{q} \leq \sqrt{2}\,p$. Moreover, for all $x,y \in X$ and $t \geq 0$ we have
\begin{equation}\label{eq:hilbert-feature-inequality}
\|F(x)\|_H = \mathfrak{q}(x), \quad \mathfrak{q}(x+y) \leq \|F(x)+F(y)\|_H, \quad F(tx) = tF(x),
\end{equation}
as well as
\begin{equation}\label{eq:triangle-defect-transfer}
\mathfrak{q}(x)+\mathfrak{q}(y)-\mathfrak{q}(x+y) \geq \frac{p(x)+p(y)-p(x+y)}{2\sqrt{2}}.
\end{equation}
\end{lemma}

\begin{proof}
First, we make use of Lemma \ref{lem:qc}(i) which gives
$$
\sum_{j=1}^\infty w_j\,q_{c_j}(x)^2 \leq \sum_{j=1}^\infty 2^{-j}p(x)^2 = p(x)^2.
$$
Hence, $F$ is well defined and has the two properties $\|F(x)\|_H = \mathfrak{q}(x)$ and $p \leq \mathfrak{q} \leq \sqrt{2}\,p$. Moreover, Lemma \ref{lem:qc}(ii) gives
\begin{align*}
\mathfrak{q}(x+y)^2 & \leq (p(x)+p(y))^2+\sum_{j=1}^\infty w_j\,\|\Phi_{c_j}(x)+\Phi_{c_j}(y)\|_{L^2(U)}^2\\
 &= \|F(x)+F(y)\|_H^2.
\end{align*}
The homogeneity of $F$ follows from Lemma \ref{lem:qc}(iii), which proves \eqref{eq:hilbert-feature-inequality}. In particular, $\mathfrak{q}(x+y) \leq \|F(x)\|_H+\|F(y)\|_H = \mathfrak{q}(x)+\mathfrak{q}(y)$. The absolute homogeneity and the positive definiteness of $\mathfrak{q}$ follow from the corresponding properties of $p$ and $q_{c_j}$. Hence, $\mathfrak{q}$ is a norm.

It remains to prove \eqref{eq:triangle-defect-transfer}. If $x=y=0$, both sides of \eqref{eq:triangle-defect-transfer} vanish. Otherwise, we have $\mathfrak{q}(x)+\mathfrak{q}(y)>0$, and the triangle inequality in $H$ and Lemma \ref{lem:qc}(ii) give
$$
\begin{aligned}
(\mathfrak{q}(x)+\mathfrak{q}(y))^2-\mathfrak{q}(x+y)^2 &\geq \|F(x)+F(y)\|_H^2-\mathfrak{q}(x+y)^2 \\
&\geq (p(x)+p(y))^2-p(x+y)^2.
\end{aligned}
$$
We divide both sides of the latter inequality by $\mathfrak{q}(x)+\mathfrak{q}(y)+\mathfrak{q}(x+y) \leq 2(\mathfrak{q}(x)+\mathfrak{q}(y))$. Since $p(x)+p(y)-p(x+y) \geq 0$ and $p(x)+p(y) \geq (\mathfrak{q}(x)+\mathfrak{q}(y))/\sqrt{2}$, this yields
$$
\begin{aligned}
\mathfrak{q}(x)+\mathfrak{q}(y)-\mathfrak{q}(x+y) &\geq \frac{(p(x)+p(y)-p(x+y))(p(x)+p(y))}{2(\mathfrak{q}(x)+\mathfrak{q}(y))} \\
&\geq \frac{p(x)+p(y)-p(x+y)}{2\sqrt{2}},
\end{aligned}
$$
and therefore \eqref{eq:triangle-defect-transfer} follows.
\end{proof}

\subsection{Proof of Theorem \texorpdfstring{\ref{thm:reflexivity-characterization}}{1.3}}\label{subsec:proof-thm13}

We are now in a position to construct, for every separable reflexive Banach space, a single equivalent norm with the $\A$-Fekete property for all $\A$ satisfying the defect condition. Note that no assumption on the dimension of $X$ is needed.

\begin{theorem}\label{thm:sufficiency}
Let $X$ be a separable reflexive Banach space over $\K$. Then, there exists a norm $\mathfrak{q}$ on $X$, equivalent to the norm of $X$, such that $\mathfrak{q}$ has the $\A$-Fekete property for every $\A \subseteq \Delta$ satisfying condition \eqref{eq:defect-condition}.
\end{theorem}

\begin{proof}
By Lemma \ref{lem:2rotund}, there exists a $2$-rotund norm $p$ on $X$ which is equivalent to the norm of $X$. Since $X$ is separable, there exists a sequence $(c_j)_{j \in \N}$ which is dense in $X$. Let $\mathfrak{q}$, $H$ and $F$ be defined as in Subsection \ref{subsec:translated-norms}. By Lemma \ref{lem:N-F}, the map $\mathfrak{q}$ is a norm which is equivalent to $p$, and hence to the norm of $X$. Notice that $\mathfrak{q}$ does not depend on $\A$.

Now let $\A \subseteq \Delta$ satisfy condition \eqref{eq:defect-condition} and let $(v_n)_{n \in \N} \subseteq X$ be $\A$-admissible with respect to $\mathfrak{q}$. We show that $(v_n/n)_{n \in \N} $ converges in norm.

\smallskip
\textbf{Step 1.} For $(n,m) \in \A$, the $\A$-admissibility of $(v_n)_{n \in \N} $ together with \eqref{eq:hilbert-feature-inequality} yields the estimate
$$
\|F(v_{n+m})\|_H = \mathfrak{q}(v_{n+m}) \leq \mathfrak{q}(v_n+v_m) \leq \|F(v_n)+F(v_m)\|_H.
$$
Hence, $(F(v_n))_{n \in \N} $ is an $\A$-admissible sequence in the real Hilbert space $H$. By Proposition \ref{prop:hilbert-sufficiency}, there exists $h = (h_0,(h_j)_{j \in \N}) \in H$ such that $F(v_n)/n \to h$ in $H$.

Let $a_n = \mathfrak{q}(v_n) = \|F(v_n)\|_H$, and let $\pi_n = p(v_n)$, which is the first coordinate of $F(v_n)$. Then, we have
\begin{equation}\label{eq:two-scalar-limits}
\frac{a_n}{n} \longrightarrow \alpha = \|h\|_H, \quad \frac{\pi_n}{n} \longrightarrow \beta = h_0, \quad \beta \leq \alpha \leq \sqrt{2}\,\beta,
\end{equation}
where the last two inequalities follow from $\pi_n \leq a_n \leq \sqrt{2}\,\pi_n$. If $\beta=0$, then $p(v_n/n) \to 0$, and hence $v_n/n \to 0$ in norm. Thus, we may assume that $\beta>0$, in which case there exists $n_0 \in \N$ such that $\pi_n>0$ for all $n \geq n_0$.

\smallskip
\textbf{Step 2.} For $n \geq n_0$, define $y_n$ via $y_n = v_n/\pi_n$. Then $p(y_n)=1$, and the homogeneity of $F$ in \eqref{eq:hilbert-feature-inequality} gives
$$
F(y_n) = \frac{F(v_n)/n}{\pi_n/n} \longrightarrow \frac{h}{\beta}, \quad n \to \infty.
$$
Fix $j\in\N$ and set $g_n(\lambda)=p(y_n+\lambda c_j)$.
The convergence of $F(y_n)$ implies that $(g_n)$ is Cauchy in
$L^2(U)$. Denoting by $B_r$ the ball $\{\lambda\in U:|\lambda-1|<r\}$
and using $|g_n(\lambda)-g_n(1)|\leq r p(c_j)$ on $B_r$ gives
$$
|g_n(1)-g_m(1)|\leq 2r p(c_j)
+\mu(B_r)^{-1/2}\|g_n-g_m\|_{L^2(U)}.
$$
Taking the $\limsup$ as $n,m\to\infty$ and then letting
$r\to0$ shows that $p(y_n+c_j)$ converges. Since the functions
$x\mapsto p(y_n+x)$ are $1$-Lipschitz and since the sequence
$(c_j)_{j \in \N} $ is dense in $X$, it follows that
\begin{equation}\label{eq:translated-type-limit}
\tau(x)=\lim_{n\to\infty}p(x+y_n)
\quad\text{exists for every }x\in X.
\end{equation}

\smallskip
\textbf{Step 3. Convergence of $(y_n)_{n \in \N} $.} Suppose by contradiction that $(y_n)_{n \in \N} $ does not converge in norm. By \eqref{eq:translated-type-limit}, we can apply Lemma \ref{lem:translated-limits} to $(y_n)_{n \geq n_0}$ and obtain the bound $\limsup_{n \to \infty} \tau(y_n)<2$. Since \eqref{eq:two-scalar-limits} gives $\alpha \geq \beta>0$, we can choose $\delta \in (0,1]$ and $K \geq n_0$ such that
\begin{equation}\label{eq:tail-conditions}
\tau(y_n) \leq 2-2\delta \quad \text{and} \quad a_n \geq \frac{\alpha n}{2}, \quad n \geq K.
\end{equation}

Now fix $n \geq K$. By \eqref{eq:translated-type-limit} and \eqref{eq:tail-conditions}, there exists $m_0(n) \geq n$ such that
\begin{equation}\label{eq:base-norm-row-gap}
p(y_n+y_m) \leq 2-\delta, \quad m \geq m_0(n).
\end{equation}
We claim that
\begin{equation}\label{eq:defect-lower}
D_a(n,m) \geq \frac{\delta\alpha}{8}\,n \quad \text{for all } m \geq m_0(n) \text{ with } (n,m) \in \A.
\end{equation}
Let $m\geq m_0(n)$ with $(n,m)\in\A$ and set
$\sigma=\min\{\pi_n,\pi_m\}$. By Maligranda's sharpened triangle
inequality \cite[Theorem~1.1, (1.1)]{maligranda2008}, it follows that
$$
\pi_n+\pi_m-p(v_n+v_m)\geq \sigma\big(2-p(y_n+y_m)\big)
\geq\delta\sigma,
$$
where for the last inequality we used \eqref{eq:base-norm-row-gap}.
Hence, by admissibility and \eqref{eq:triangle-defect-transfer},
\begin{align*}
D_a(n,m) &\geq \frac{\pi_n+\pi_m-p(v_n+v_m)}{2\sqrt{2}}\\
  &\geq \frac{\delta\sigma}{2\sqrt{2}}\\
  &\geq \frac{\delta}{4}\min\{a_n,a_m\}\\
  &\geq \frac{\delta\alpha}{8}\,n.
\end{align*}
This proves \eqref{eq:defect-lower}.

By \eqref{eq:defect-lower}, for every $n \geq K$ the set
\begin{equation}
\{m \in \N : (n,m) \in \A, \ D_a(n,m)<(\delta\alpha/8)\,n\},
\end{equation}
is contained in $\{1,\dots,m_0(n)-1\}$, and hence it is finite. On the other hand, the sequence $a=(a_n)_{n \in \N} $ is $\A$-subadditive. Moreover, by \eqref{eq:two-scalar-limits}, there exists $C>0$ with $a_n \leq Cn$ for all $n \in \N$. This contradicts Lemma \ref{lem:small-defect-row} with $\varepsilon = \delta\alpha/8$. Hence, $(y_n)_{n \in \N} $ converges in norm to some $y \in X$.

\smallskip
\textbf{Step 4. Conclusion.} It follows from Steps 1 and 3 that
$$
\frac{v_n}{n} = \frac{\pi_n}{n}\,y_n \longrightarrow \beta y, \quad n \to \infty.
$$
Hence, $(v_n/n)_{n \in \N} $ converges in norm. Since $\A$ was arbitrary and $\mathfrak{q}$ does not depend on $\A$, the proof is complete.
\end{proof}

Theorem \ref{thm:sufficiency}, combined with Corollary \ref{cor:plane-obstruction} and Lemma \ref{lem:james-obstruction}, yields Theorem \ref{thm:reflexivity-characterization}.

\begin{proof}[Proof of Theorem \ref{thm:reflexivity-characterization}]
Let $X$ be a separable Banach space with $\dim X \geq 2$ and let $\A \subseteq \Delta$.

We first prove (i)$\Rightarrow$(ii). Let $q$ be a norm on $X$ which is equivalent to the norm of $X$ and has the $\A$-Fekete property. Since $\dim_\R X \geq \dim X \geq 2$, Corollary \ref{cor:plane-obstruction}, applied to the normed space $(X,q)$, shows that $\A$ satisfies condition \eqref{eq:defect-condition}. Moreover, $X$ is reflexive (since otherwise $q$ would not have the $\A$-Fekete property by Lemma \ref{lem:james-obstruction}). The implication (ii)$\Rightarrow$(i) and the additional assertion of the theorem, namely that there is one equivalent norm on $X$ which has the $\A$-Fekete property for every $\A \subseteq \Delta$ satisfying condition \eqref{eq:defect-condition}, follow from Theorem \ref{thm:sufficiency}.
\end{proof}

\section{Proof of Theorem \texorpdfstring{\ref{thm:transfer}}{1.1} and further consequences}\label{sec:consequences}

In this section, we prove Theorem \ref{thm:transfer} and derive further consequences of our results. In particular, we show that if $\dim_\R X \geq 2$, then the inequalities \eqref{eq:fekete} restricted to $n \leq m \leq g(n)$ never imply the convergence of $(v_n/n)_{n \in \N} $, regardless of the choice of $g : \N \to \N$.

\subsection{Proof of Theorem \texorpdfstring{\ref{thm:transfer}}{1.1}}

Theorem \ref{thm:transfer} follows by combining Corollary \ref{cor:plane-obstruction} with the results of Section \ref{sec:renorming}.

\begin{proof}[Proof of Theorem \ref{thm:transfer}]
We first prove (i)$\Rightarrow$(ii). To do so, let $q$ be a norm on $X$ which is equivalent to the norm of $X$ and has the $\A$-Fekete property. Every sequence $(v_n)_{n \in \N} $ with $q(v_{n+m}) \leq q(v_n+v_m)$ for all $(n,m) \in \Delta$ is $\A$-admissible with respect to $q$. Hence, $q$ has the unrestricted Fekete property. Moreover, since $\dim_\R X \geq \dim X \geq 2$, Corollary \ref{cor:plane-obstruction}, applied to the normed space $(X,q)$, shows that $\A$ satisfies condition \eqref{eq:defect-condition}.

We now prove (ii)$\Rightarrow$(i). Let $q$ be a norm on $X$ which is equivalent to the norm of $X$ and has the unrestricted Fekete property. By Lemma \ref{lem:james-obstruction}, the space $X$ is reflexive. Since $\A$ satisfies condition \eqref{eq:defect-condition}, Theorem \ref{thm:sufficiency} provides a norm on $X$ which is equivalent to the norm of $X$ and has the $\A$-Fekete property.
\end{proof}

\subsection{Restrictions with finite rows}

Recall from Subsection \ref{subsec:rows} that the $n$-th row of $\A$ is given by the set $\{m \in \N : (n,m) \in \A\}$. Now, suppose that $\A$ satisfies condition \eqref{eq:defect-condition}. Applying Lemma \ref{lem:small-defect-row} to the zero sequence and $\varepsilon=1$ shows that for every $K \in \N$ there exists $n \geq K$ such that the $n$-th row of $\A$ is infinite. Thus, infinitely many rows of $\A$ are infinite.

\begin{corollary}\label{cor:finite-rows}
Let $\A \subseteq \Delta$ be such that only finitely many rows of $\A$ are infinite. Then $\A$ does not satisfy condition \eqref{eq:defect-condition}. Moreover, for every normed space $X$ with $\dim_\R X \geq 2$ there exists an $\A$-admissible sequence $(v_n)_{n \in \N}  \subseteq X$ with the following properties:
\begin{enumerate}[\rm (i)]
\item The sequence $(\|v_n\|/n)_{n \in \N} $ is strictly decreasing and converges to $1$.
\item The sequence $(\|v_n\|)_{n \in \N} $ is subadditive, i.e., $\|v_{n+m}\| \leq \|v_n\|+\|v_m\|$ for all $n,m \in \N$.
\item The sequence $(v_n/n)_{n \in \N} $ has two distinct accumulation points.
\end{enumerate}
\end{corollary}

\begin{proof}
Choose $K\in\N$ such that the $n$-th row of $\A$ is finite for
every $n\geq K$. Applying Lemma~\ref{lem:block-construction} with
$a=0$ and $\varepsilon=1$ gives a sequence
$b_n=(n/\kappa(n)+1)/3$ in $\Pro(\A)$ satisfying \eqref{eq:failure}.
In particular, $\A$ fails the defect condition.

Apply Proposition~\ref{prop:plane-obstruction} in a two-dimensional
real subspace of $X$. The resulting sequence $(v_n)_{n \in \N} $ has
$a_n:=\|v_n\|=n+b_n$, and $(v_n/a_n)_{n \in \N} $ has two distinct accumulation
points. Since $\kappa$ is nondecreasing and tends to infinity,
$$
\frac{a_n}{n}=1+\frac{1}{3}\left(\frac{1}{\kappa(n)}+\frac{1}{n}\right)
$$
satisfies $a_{n+1}/(n+1)<a_n/n$ for all $n\in\N$ and $\lim_{n\to\infty}a_n/n=1$. This proves (i).

Comparing $a_{n+m}/(n+m)$
with $a_n/n$ and $a_m/m$ gives $a_{n+m}\leq a_n+a_m$ for all
$n,m\in\N$, which proves (ii).

Finally,
$$
\left\|\frac{v_n}{n}-\frac{v_n}{a_n}\right\|
=\left|\frac{a_n}{n}-1\right|\to0,
$$
which shows that the two sequences have the same accumulation points.
\end{proof}

As a special case, we obtain that the inequalities \eqref{eq:fekete} restricted to $n \leq m \leq g(n)$ never force convergence, regardless of the choice of $g$.

\begin{corollary}\label{cor:upper-bounds}
Let $X$ be a normed space with $\dim_\R X \geq 2$, $g : \N \to \N$ and $I \subseteq \N$ be finite. Define
$$
\A = \{(n,m) \in \Delta : m \leq g(n)\} \cup \{(n,m) \in \Delta : n \in I\}.
$$
Then, there exists an $\A$-admissible sequence $(v_n)_{n \in \N}  \subseteq X$ such that $(v_n/n)_{n \in \N} $ does not converge. In particular, there exists no function $g : \N \to \N$ such that the inequalities
$$
\|v_{n+m}\| \leq \|v_n+v_m\|, \quad n \leq m \leq g(n),
$$
imply the convergence of $(v_n/n)_{n \in \N} $ for every sequence $(v_n)_{n \in \N}  \subseteq X$.
\end{corollary}

\begin{proof}
For $n \notin I$, the $n$-th row of $\A$ is contained in $\{n,\dots,g(n)\}$, and hence it is finite. Thus, only finitely many rows of $\A$ are infinite, and the claim follows from Corollary \ref{cor:finite-rows}. The last assertion is the case $I=\emptyset$.
\end{proof}

Corollary \ref{cor:upper-bounds} answers a question of Kulikov and Shao \cite[Question 3.2]{kulikov2024feketeslemmabanachspaces} in the negative for every normed space $X$ with $\dim_\R X \geq 2$, in particular for every uniformly convex Banach space. Moreover, Corollary \ref{cor:finite-rows} shows that the convergence of $(\|v_n\|/n)_{n \in \N} $ does not imply the convergence of $(v_n/n)_{n \in \N} $, even if $(\|v_n\|)_{n \in \N} $ is subadditive.

\section{AI disclosure}

The proof of Theorem \ref{thm:sufficiency}, which gives the implication (ii)$\Rightarrow$(i) in Theorem \ref{thm:reflexivity-characterization}, was obtained with the assistance of the large language models GPT-5.5 and GPT-5.6. To that end, we suggested to the models to base the proof on a rotundity property of a suitable equivalent norm and supplied them with the paper of Odell and Schlumprecht \cite{odell1998asymptotic}. A proof emerged from a sustained conversation with both models, during which new ideas arose and the picture of the problem gradually became complete. The original proof was lengthy, and we subsequently simplified it significantly.

\section{Acknowledgments}
V.B. acknowledges support from VILLUM FONDEN via the QMATH Centre of Excellence (Grant No.~10059), and the Novo Nordisk Foundation (Grant No.~NNF20OC0059939, “Quantum for Life”).
L.L.~is grateful to the Azrieli Foundation for the award of an Azrieli Fellowship and acknowledges the support of this research by ISF Grant No.~854/25.
S.G.~is supported by Grant 334466 of the Research Council of Norway.


\begin{thebibliography}{9}

\bibitem{deBruijnErdos1952}
N.~G. de Bruijn and P.~Erd\H{o}s,
\emph{Some linear and some quadratic recursion formulas.~II},
Indag. Math. \textbf{14} (1952), 152--163.

\bibitem{fekete1923}
M.~Fekete,
\emph{\"Uber die Verteilung der Wurzeln bei gewissen algebraischen
Gleichungen mit ganzzahligen Koeffizienten},
Math. Z. \textbf{17} (1923), 228--249.

\bibitem{furedi2020nearly}
Z. F\"uredi and I. Ruzsa, \emph{Nearly subadditive sequences},
Acta Math.\ Hung. \textbf{161} (2020), 401--411.

\bibitem{james1972reflexivity}
R.~C. James,
\emph{Reflexivity and the sup of linear functionals},
Israel J. Math. \textbf{13} (1972), 289--300.
\href{https://doi.org/10.1007/BF02762803}{doi:10.1007/BF02762803}.

\bibitem{kulikov2024feketeslemmabanachspaces}
A.~Kulikov and F.~Shao,
\emph{Fekete's lemma in Banach spaces},
\href{https://arxiv.org/abs/2411.17380}{arXiv:2411.17380} (2024).

\bibitem{odell1998asymptotic}
E.~Odell and Th.~Schlumprecht,
\emph{Asymptotic properties of Banach spaces under renormings},
J. Amer. Math. Soc. \textbf{11} (1998), no.~1, 175--188.
\href{https://doi.org/10.1090/S0894-0347-98-00251-3}
{doi:10.1090/S0894-0347-98-00251-3}.

\bibitem{brezis2011}
H.~Brezis,
\emph{Functional Analysis, Sobolev Spaces and Partial Differential Equations},
Universitext, Springer, New York, 2011.
\href{https://doi.org/10.1007/978-0-387-70914-7}{doi:10.1007/978-0-387-70914-7}.

\bibitem{dilworth2024}
S.~Dilworth, D.~Kutzarova and P.~Motakis,
\emph{$2$-Rotund norms for unconditional and symmetric sequence spaces},
Banach J. Math. Anal. \textbf{18} (2024), article 74.
\href{https://doi.org/10.1007/s43037-024-00379-1}{doi:10.1007/s43037-024-00379-1}.

\bibitem{kunduJain2006}
S.~Kundu and T.~Jain,
\emph{Atsuji spaces: equivalent conditions},
Topology Proc. \textbf{30} (2006), no.~1, 301--325.

\bibitem{maligranda2008}
L.~Maligranda,
\emph{Some remarks on the triangle inequality for norms},
Banach J. Math. Anal. \textbf{2} (2008), no.~2, 31--41.
\href{https://doi.org/10.15352/bjma/1240336290}{doi:10.15352/bjma/1240336290}.

\end{thebibliography}
\end{document}